\documentclass[11pt]{article}
\usepackage{epsfig}

\def\be{\begin{equation}}
\def\ee{\end{equation}}
\def\bea{\begin{eqnarray}}
\def\eea{\end{eqnarray}}
\def\bes{\begin{eqnarray*}}
\def\ees{\end{eqnarray*}}

\def\<{\langle}
\def\>{\rangle}
\def\lb{\label}
\def\bs{\setminus}
\def\pt{\partial}
\def\grad{{\rm grad}}
\def\R{{\bf R}}
\def\C{{\bf C}}
\def\Z{{\bf Z}}
\def\N{{\bf N}}
\def\U{{\bf U}}
\def\Q{{\bf Q}}

\def\CP{{\bf CP}}

\def\ga{{\gamma}}

\def\ka{{\kappa}}
\def\th{{\theta}}

\def\lm{{\lambda}}
\def\Lm{{\Lambda}}

\def\im{{\rm im}}

\def\rank{{\rm rank}}
\def\Sp{{\rm Sp}}

\def\hb{\vrule height0.18cm width0.14cm $\,$}

\title{ Existence of closed geodesics on Finsler $n$-spheres}

\author{Wei Wang\thanks{Partially supported by National Natural
Science Foundation of China No.10801002, China Postdoctoral Science Foundation
No.20070420264 and LMAM in Peking University.
E-mail: alexanderweiwang@yahoo.com.cn, wangwei@math.pku.edu.cn  }\\
School of Mathematical Science \\ Peking University, Beijing 100871 \\
PEOPLES REPUBLIC OF CHINA \\ }

\date{April 30th, 2009}
\begin{document}

\maketitle

\begin{abstract}
{\it In this paper, we prove that on every Finsler  $n$-sphere
$(S^n,\,F)$ with reversibility $\lambda$
satisfying $F^2<(\frac{\lambda+1}{\lambda})^2g_0$
and $l(S^n, F)\ge \pi(1+\frac{1}{\lambda})$,
there always exist at least $n$ prime closed geodesics without self-intersections,
where $g_0$ is the standard Riemannian metric on $S^n$
with constant curvature $1$ and $l(S^n, F)$ is the length
of a shortest geodesic loop on $(S^n, F)$.
We also study the stability of these closed geodesics. }
\end{abstract}



{\bf Running head}: Closed geodesics on Finsler spheres

\renewcommand{\theequation}{\thesection.\arabic{equation}}
\renewcommand{\thefigure}{\thesection.\arabic{figure}}

\setcounter{equation}{0}%\setcounter{figure}{0}
\section{Introduction and main results}%{Section 1}

There is a famous conjecture in Riemannian geometry which
claims there exist infinitely many prime
closed geodesics on any Riemannian manifold. This conjecture
has been proved except for CROSS's (compact rank one symmetric
spaces). The results of J. Franks \cite{Fra} in 1992 and V. Bangert \cite{Ban}
in 1993 imply this conjecture is true for any Riemannian 2-sphere.
But once one move to the Finsler case, the conjecture
becomes false. It was quite surprising when A. Katok \cite{Kat}  in 1973 found
some non-symmetric Finsler metrics on CROSS's with only finitely
many prime closed geodesics and all closed geodesics are
non-degenerate and elliptic. In Katok's examples the spheres $S^{2n}$ and
$S^{2n-1}$ have precisely $2n$ closed geodesics.
Based on Katok's work, D. V. Anosov \cite{Ano} in 1974's ICM
conjectured that there exist at least $2[\frac{n+1}{2}]$ prime
closed geodesics on any Finsler  $n$-sphere $(S^n,\,F)$.
In \cite{Zil}, W. Ziller conjectured the number of
prime closed geodesics on any Finsler  $n$-sphere $(S^n,\,F)$
is at least $n$. This paper is devoted to a study on the number
prime closed geodesics on Finsler $n$-spheres.
Let us recall firstly the definition of the Finsler
metrics.

{\bf Definition 1.1.} (cf. \cite{She}) {\it Let $M$ be a finite
dimensional smooth manifold. A function $F:TM\to [0,+\infty)$ is a {\rm
Finsler metric} if it satisfies

(F1) $F$ is $C^{\infty}$ on $TM\bs\{0\}$,

(F2) $F(x,\lm y) = \lm F(x,y)$ for all $x\in M$, $y\in T_xM$ and
$\lm>0$,

(F3) For every $y\in T_xM\bs\{0\}$, the quadratic form
$$ g_{x,y}(u,v) \equiv
         \frac{1}{2}\frac{\pt^2}{\pt s\pt t}F^2(x,y+su+tv)|_{t=s=0},
         \qquad \forall u, v\in T_xM, $$
is positive definite.

In this case, $(M,F)$ is called a {\rm Finsler manifold}. $F$ is
{\rm symmetric} if $F(x,-y)=F(x,y)$ holds for all $x\in M$ and
$y\in T_xM$. $F$ is {\rm Riemannian} if $F(x,y)^2=\frac{1}{2}G(x)y\cdot
y$ for some symmetric positive definite matrix function $G(x)\in
GL(T_xM)$ depending on $x\in M$ smoothly. }

A closed curve in a Finsler manifold is a closed geodesic if it is
locally the shortest path connecting any two nearby points on this
curve (cf. \cite{She}). As usual, on any Finsler manifold
$(M, F)$, a closed geodesic $c:S^1=\R/\Z\to M$ is {\it prime}
if it is not a multiple covering (i.e., iteration) of any other
closed geodesics. Here the $m$-th iteration $c^m$ of $c$ is defined
by $c^m(t)=c(mt)$. The inverse curve $c^{-1}$ of $c$ is defined by
$c^{-1}(t)=c(1-t)$ for $t\in \R$.  Note that unlike Riemannian manifold,
the inverse curve $c^{-1}$ of a closed geodesic $c$
on a non-symmetric Finsler manifold need not be a geodesic.
We call two prime closed geodesics
$c$ and $d$ {\it distinct} if there is no $\th\in (0,1)$ such that
$c(t)=d(t+\th)$ for all $t\in\R$.
We shall omit the word {\it distinct} when we talk about more than one prime closed geodesic.
On a symmetric Finsler (or Riemannian) manifold, two closed geodesics
$c$ and $d$ are called { \it geometrically distinct} if $
c(S^1)\neq d(S^1)$, i.e., their image sets in $M$ are distinct.

For a closed geodesic $c$ on $(M,\,F)$, denote by $P_c$
the linearized Poincar\'{e} map of $c$  (cf. p.143 of \cite{Zil}).
Then $P_c\in \Sp(2n-2)$ is a symplectic matrix.
For any $M\in \Sp(2k)$, we define the {\it elliptic height } $e(M)$
of $M$ to be the total algebraic multiplicity of all eigenvalues of
$M$ on the unit circle $\U=\{z\in\C|\; |z|=1\}$ in the complex plane
$\C$. Since $M$ is symplectic, $e(M)$ is even and $0\le e(M)\le 2k$.
Then $c$ is called {\it hyperbolic} if all the eigenvalues of $P_c$ avoid
the unit circle in $\C$, i.e., $e(P_c)=0$; {\it elliptic} if all the eigenvalues of
$P_c$ are on the unit circle, i.e., $e(P_c)=2(n-1)$.

Following H-B. Rademacher in
\cite{Rad4}, the reversibility $\lambda=\lambda(M,\,F)$ of a compact
Finsler manifold $(M,\,F)$ is defied to be
$$\lambda:=\max\{F(-X)\,|\,X\in TM, \,F(X)=1\}\ge 1.$$

We are aware of a number of results concerning closed geodesics
on spheres. According to the classical theorem of Lyusternik-Fet \cite{LyF} from 1951, there exists at least one
closed geodesic on every compact Riemannian manifold. The proof of this theorem
is variational and carries over to the Finsler case. Motivated by the
work \cite{Kli1} of W. Klingenberg in 1969, W. Ballmann, G. Thorbergsson and
W. Ziller studied in \cite{BTZ1} and \cite{BTZ2} of 1982-83 the existence and
stability of closed geodesics on positively curved CROSS's.
In \cite{Hin} of 1984, N. Hingston proved that a Riemannian metric
on a sphere all of whose closed geodesics are hyperbolic carries infinitely
many geometrically distinct closed geodesics.
By the results of J. Franks in \cite{Fra} of 1992  and V. Bangert
in \cite{Ban} of 1993, there are infinitely many geometrically
distinct closed geodesics for any Riemannian metric on $S^2$.
In \cite{Rad5}, H.-B. Rademacher studied the existence and
stability of closed geodesics on positively curved Finsler manifolds.
In \cite{BaL}, V. Bangert and Y. Long proved that on any Finsler 2-sphere $(S^2, F)$,
there exist at least two prime closed geodesics.
In \cite{LoW} of Y. Long and the author, they further
proved the existence of at least two irrationally elliptic
prime closed geodesics on every Finsler  $2$-sphere $(S^2,\,F)$
provided the number of prime closed geodesics is finite.
In \cite{Wang}, the author proved there exist three prime
closed geodesics on any $(S^3, F)$ satisfying
$(\frac{\lambda}{\lambda+1})^2<K\le 1$,
where $K$ is the flag curvature of $(S^3, F)$.

The following are the main results in this paper:

{\bf Theorem 1.2.} {\it On every Finsler  $n$-sphere
$(S^n,\,F)$ with reversibility $\lambda$
satisfying $F^2<(\frac{\lambda+1}{\lambda})^2g_0$
and $l(S^n, F)\ge \pi(1+\frac{1}{\lambda})$,
there always exist at least $n$ prime closed geodesics without self-intersections,
where $g_0$ is the standard Riemannian metric on $S^n$
with constant curvature $1$ and $l(S^n, F)$ is the length
of a shortest geodesic loop on $(S^n, F)$. }

{\bf Theorem 1.3.} {\it On every Finsler  $n$-sphere
$(S^n,\,F)$ with reversibility $\lambda$
satisfying $F^2<(\frac{\lambda+1}{\lambda})^2g_0$
and $\left(\frac{\lambda}{\lambda+1}\right)^2<K\le 1$
($0<K\le 1$ if n is even), there always exist at least
$n$ prime closed geodesics without self-intersections,
where $g_0$ is given by Theorem 1.2. }

{\bf Theorem 1.4.} {\it On every Finsler  $n$-sphere
$(S^n,\,F)$ with reversibility $\lambda$
satisfying $F^2<(\frac{2(n-1)}{2n-1})^2g_0$
and $(\frac{2n-3}{n-1}\frac{\lambda}{\lambda+1}
)^2<K\le 1$,
there exist at least $2[\frac{n}{2}]$ non-hyperbolic prime closed geodesics,
where $g_0$ is given by Theorem 1.2. }

{\bf Remark 1.5.}
The proof of these theorems is motivated by \cite{EL},
\cite{Rad3}, \cite{BTZ1} and \cite{BTZ2}. We use the
$S^1$-equivariant Morse theory
to obtain $n$ critical values of the energy functional
$E$ on the space pair $(\Lambda,\,\Lambda^0)$,
where $\Lambda$ is the free loop space of $S^n$ and $\Lambda^0$
is its subspace consisting of constant point curves.
Then we use the assumptions to show that these critical values correspond
to $n$ distinct prime closed geodesics. The stability results are obtained
by considering the indices of iterations of these closed geodesics,
cf. \cite{Bot2}.
Note that the methods in \cite{BTZ1} and \cite{BTZ2} can't be used for
non-symmetric Finsler metrics because the lack of $\Z_2$-symmetry
for $E$.

In this paper, let $\N$, $\N_0$, $\Z$, $\Q$, $\R$, and $\C$ denote
the sets of natural integers, non-negative integers, integers,
rational numbers, real numbers, and complex numbers respectively.
We use only singular homology modules with $\Q$-coefficients.
We denote by $[a]=\max\{k\in\Z\,|\,k\le a\}$ for any $a\in\R$.

\setcounter{equation}{0}%\setcounter{figure}{0}
\section{Critical point theory for closed geodesics}%{Section 1}

In this section, we will study critical point theory for closed geodesics.

On a compact Finsler manifold $(M,F)$, we choose an auxiliary Riemannian
metric. This endows the space $\Lambda=\Lambda M$ of $H^1$-maps
$\gamma:S^1\rightarrow M$ with a natural Riemannian Hilbert manifold structure
on which the group $S^1=\R/\Z$ acts continuously
by isometries, cf. \cite{Kli2}, Chapters 1 and 2. This action is
defined by translating the parameter, i.e.,
$$ (s\cdot\gamma)(t)=\gamma(t+s) \qquad $$
for all $\gamma\in\Lm$ and $s,t\in S^1$.
The Finsler metric $F$ defines an energy functional $E$ and a length
functional $L$ on $\Lambda$ by
\be E(\gamma)=\frac{1}{2}\int_{S^1}F(\dot{\gamma}(t))^2dt,
 \quad L(\gamma) = \int_{S^1}F(\dot{\gamma}(t))dt.  \lb{2.1}\ee
Both functionals are invariant under the $S^1$-action.
By \cite{Mer}, the functional $E$ is $C^{1, 1}$ on $\Lm$ and
satisfies the Palais-Smale condition. Thus we can apply the
$S^1$-deformation theorems in \cite{Cha} and \cite{MaW}.
The critical points
of $E$ of positive energies are precisely the closed geodesics $c:S^1\to M$
of the Finsler structure. If $c\in\Lambda$ is a closed geodesic then $c$ is
a regular curve, i.e., $\dot{c}(t)\not= 0$ for all $t\in S^1$, and this implies
that the formal second differential $E''(c)$ of $E$ at $c$ exists.
As usual we define the index $i(c)$ of $c$ as the maximal dimension of
subspaces of $T_c \Lambda$ on which $E^{\prime\prime}(c)$ is negative definite, and the
nullity $\nu(c)$ of $c$ so that $\nu(c)+1$ is the dimension of the null
space of $E^{\prime\prime}(c)$. In fact, one can use a finite dimension
approximation $\Delta$ of $\Lambda$ as in \cite{Rad2} such that
$F|_\Delta$ is $C^\infty$ and define the index and nullity,
one can prove those equal to the above defined ones.

For $m\in\N$ we denote the $m$-fold iteration map
$\phi^m:\Lambda\rightarrow\Lambda$ by \be \phi^m(\ga)(t)=\ga(mt)
\qquad \forall\,\ga\in\Lm, t\in S^1. \lb{2.2}\ee We also use the
notation $\phi^m(\gamma)=\gamma^m$. For a closed geodesic $c$, the
average index is defined by \be
\hat{i}(c)=\lim_{m\rightarrow\infty}\frac{i(c^m)}{m}. \lb{2.3}\ee

If $\gamma\in\Lambda$ is not constant then the multiplicity
$m(\gamma)$ of $\gamma$ is the order of the isotropy group $\{s\in
S^1\mid s\cdot\gamma=\gamma\}$. If $m(\gamma)=1$ then $\gamma$ is
called {\it prime}. Hence $m(\gamma)=m$ if and only if there exists a
prime curve $\tilde{\gamma}\in\Lambda$ such that
$\gamma=\tilde{\gamma}^m$.

In this paper for $\ka\in \R$ we denote by
\be \Lm^{\ka}=\{d\in \Lm\,|\,E(d)\le \ka\},    \lb{2.4}\ee
For a closed geodesic $c$ we set
$$ \Lm(c)=\{\ga\in\Lm\mid E(\ga)\le E(c)\}. $$

We have the following property for the space pair
$H_{S^1, \ast}(\Lambda S^n, \Lambda^0 S^n)$,
where $H_{S^1, \ast}$ is the $S^1$-equivariant homology
in the sense of A. Borel, cf. \cite{Bor}.

{\bf Proposition 2.1.} {\it There are $n$ subordinate
nonzero homology classes $\sigma_k\in H_{S^1, n-1+2k}(\Lambda S^n, \Lambda^0 S^n)$ for
$0\le k\le n-1$.}

{\bf Proof.} Let $g_0$ be the standard metric on $S^n$ with constant
curvature $1$ and $E_0$ be the corresponding energy functional. Then
by \S4 of \cite{Hin}, $H_{S^1, \ast}(\Lambda S^n, \Lambda^0 S^n)$ is
generated by the local critical groups \be H_{S^1,
\ast}(\Lambda_0^{2m^2\pi^2+\epsilon} S^n,
\Lambda_0^{2m^2\pi^2-\epsilon} S^n) \simeq H_{S^1,
\ast-(2m-1)(n-1)}(T_1S^n),\quad m\in\N\lb{2.5}\ee where
$\Lm_0^{\ka}=\{d\in \Lm\,|\,E_0(d)\le \ka\}$,
$T_1S^n=\{v\in TS^n| |v|_{g_0}=1\}$ and $\epsilon>0$ is sufficiently
small. In fact, this follows since the Morse series is lacunary,
thus the energy functional $E_0$ is perfect. Now we consider
$H_{S^1, \ast}(T_1S^n)$. Note that the isotropy group of the $S^1$-action
on $T_1S^n$ is $\Z_m$ for $m\in\N$. Thus it follows from Lemma 6.11
of \cite{FaR} and the universal coefficient theorem that \be H_{S^1,
\ast}(T_1S^n)\simeq H_{\ast}(T_1S^n/S^1). \lb{2.6}\ee Now consider
the case $m=1$. In this case we have an $S^1$-fibration
$S^1\rightarrow T_1S^n\rightarrow T_1S^n/S^1$. Since the
$S^1$-action is free, $T_1S^n/S^1$ is a smooth manifold. Let $e\in
H^2(T_1S^n/S^1)$ be the Euler class of the $S^1$-fibration
$T_1S^n\rightarrow T_1S^n/S^1$. Then it follows from \cite{Wei} that
$e^{n-1}\neq0$ is a generator of  $H^{2n-2}(T_1S^n/S^1)$. Denote by
$[C]$ the fundamental class of $T_1S^n/S^1$. Then we have $n$
subordinate nonzero homology classes $\alpha_k\in H_{S^1,
2k}(T_1S^n/S^1)$ for $0\le k\le n-1$ defined by $\alpha_k=[C]\cap
e^{n-1-k}$. Now we construct $\sigma_k$ from $\alpha_k$ as the
following. First note that $T_1S^n/S^1$ is a non-degenerate critical
$S^1$-manifold in the sense of \cite{Bot1}. Thus it follows from the
handle-bundle theorem in \cite{Was}, $\Lambda_0^{2\pi^2+\epsilon}
S^n$ is $S^1$-homotopic to $\Lambda_0^{2\pi^2-\epsilon} S^n$ with
the handle-bundle $DN^-$ attached along $SN^-$, where $DN^-$ is the
closed disk bundle of the negative bundle over $T_1S^n/S^1$ and
$SN^-=\partial DN^-$. In particular, $\rank DN^-=n-1$. Now we have
\be H_{S^1, \ast}(\Lambda_0^{2\pi^2+\epsilon} S^n,
\Lambda_0^{2\pi^2-\epsilon} S^n) \simeq H_{S^1, \ast}(DN^-,
SN^-)\simeq H_{S^1, \ast-(n-1)}(T_1S^n).\lb{2.7}\ee The latter
isomorphism is given by the Thom isomorphism $\Phi$. Let $f: (\Lambda
S^n)_{S^1}\equiv\Lambda S^n\times_{S^1} S^\infty\rightarrow
\CP^\infty$ be a classifying map and $\eta\in H^2(\CP^\infty)$ be
the universal first rational Chern class. Let
$\sigma_k=\Phi^{-1}(\alpha_k)=\Phi^{-1}([C]\cap e^{n-1-k})\neq0$.
Then clearly we have $\sigma_k=\Phi^{-1}([C])\cap
(f^\ast\eta)^{n-1-k}$. In fact, this follows since
the Euler class $e\in H^2(T_1S^n)$ coincide with the first
Chern class $c_1\in H^2(T_1S^n)$ of the $S^1$-bundle
$T_1S^n\rightarrow T_1S^n/S^1$; the $S^1$-action
on $T_1S^n$ is free, thus $(T_1S^n)_{S^1}$ is
homotopic to $T_1S^n/S^1$; and $DN^-$ is $S^1$-homotopic to
$T_1S^n$.

The proof of the proposition is complete.\hfill\hb

Now suppose $\sigma_k\in H_{S^1, n-1+2k}(\Lambda S^n, \Lambda^0 S^n)$ for
$0\le k\le n-1$ are given by Proposition 2.1. Let
$\lambda_k=\inf_{\gamma\in\sigma_k}\sup_{d\in\im\gamma}E(d)$.
Then we have the following

{\bf Proposition 2.2.} {\it Each $\lambda_k$ is a critical value of
$E$ and $0<\lambda_0\le\lambda_1\le\cdots\le\lambda_{n-1}$. In particular,
if $\lambda_k=\lambda_{k+1}$ for some $0\le k<n-1$, then there are infinitely
many prime closed geodesics on $(S^n, F)$.
 }

{\bf Proof.} By \cite{Mer}, the functional $E$ is $C^{1, 1}$ on $\Lm$ and
satisfies the Palais-Smale condition. Thus the negative gradient flow
of $E$ exists and is $S^1$-equivariant. Suppose $\lambda_k$ is
not a critical value of $E$, then by the Palais-Smale condition,
there exist some constants $\delta, \rho>0$ such that
$\|\grad E\|\ge\delta$ for any
$d\in E^{-1}[\lambda_k-\rho, \lambda_k+\rho]$.
Now by definition of $\lambda_k$, we have a chain
$\gamma\in\sigma_k$ such that $\sup_{d\in\im\gamma}E(d)<\lambda_k+\rho$.
Thus we can push $\gamma$ by the negative gradient flow
of $E$ down the level $\lambda_k-\rho$ to obtain $\gamma^\prime$.
Since the flow is $S^1$-equivariant, we have $[\gamma]=[\gamma^\prime]$.
This contradicts to the definition of $\lambda_k$.

Since any nonconstant geodesic loop $c$ on $(S^n , F)$ satisfies $L(c)>\varrho$ for
some $\varrho>0$. We have $E(c)\ge \frac{1}{2}L(c)^2>\frac{1}{2}\varrho^2>0$
for any nonconstant closed geodesic $c$ on $(S^n, F)$.
We claim that $\lambda_0\ge\frac{1}{2}\varrho^2$. In fact,
suppose $\lambda_0<\frac{1}{2}\varrho^2$, then
by definition of $\lambda_0$, we have a chain
$\gamma\in\sigma_0$ such that $\sup_{d\in\im\gamma}E(d)<\frac{1}{2}\varrho^2$.
Since there is no critical value of $E$ in the interval $(0, \frac{1}{2}\varrho^2)$,
we have $\sigma_0\in H_{S^1, n-1}(\Lambda^{\frac{1}{2}\varrho^2-\epsilon} S^n, \Lambda^0 S^n)=0$,
where $\epsilon>0$ is sufficiently small.
This contradicts to $\sigma_0\neq 0$.

The relation $\lambda_0\le\lambda_1\le\cdots\le\lambda_{n-1}$
is obvious by the definition of cap product.

Now suppose $\lambda_k=\lambda_{k+1}$ for some $0\le k<n-1$.
Denote by $crit(E)$ the set of critical points of $E$.
Suppose $\mathcal{U}$ is any $S^1$-invariant neighborhood
of $E^{-1}(\lambda_k)\cap crit(E)$. Then by Theorem 1.7.1 of
\cite{Cha}, there exist $\epsilon>0$ and an $S^1$-equivariant
deformation $F:[0, 1]\times\Lambda^{\lambda_k+\epsilon}\rightarrow \Lambda^{\lambda_k+\epsilon}$ such that
$F(1, \Lambda^{\lambda_k+\epsilon}\setminus\mathcal{U})\subset\Lambda^{\lambda_k-\epsilon}$.
In particular, $F$ induces a deformation
$$\widetilde{F}:[0, 1]\times(\Lambda^{\lambda_k+\epsilon})_{S^1}\rightarrow (\Lambda^{\lambda_k+\epsilon})_{S^1},
\quad
\widetilde{F}(1, (\Lambda^{\lambda_k+\epsilon}\setminus\mathcal{U})_{S^1})\subset(\Lambda^{\lambda_k-\epsilon})_{S^1}.$$
By definition of $\lambda_{k+1}$, we have a chain
$\gamma\in\sigma_{k+1}$ such that $\sup_{d\in\im\gamma}E(d)<\lambda_{k+1}+\epsilon=\lambda_k+\epsilon$.
Choose a subdivision $\gamma^{\prime}$ of $\gamma$ such that any
simplex of $\gamma^{\prime}$ is either contained in $\mathcal{U}_{S^1}$
or in $(\Lambda^{\lambda_k+\epsilon}\setminus\mathcal{U})_{S^1}$. Note that by Proposition 2.1,
we have $\sigma_k=[\gamma^{\prime}]\cap\eta$ for some
$\eta\in H^2_{S^1}(\Lambda S^n)$.
Let $i: \mathcal{U}_{S^1}\rightarrow (\Lambda S^n)_{S^1}$
and $j: (\Lambda S^n)_{S^1}\rightarrow ((\Lambda S^n)_{S^1}, \mathcal{U}_{S^1})$
be the inclusions. If $i^\ast(\eta)=0$, then we have
$\eta=j^\ast(\theta)$ for some $\theta\in H^2_{S^1}(\Lambda S^n, \mathcal{U})$
by the exactness of the cohomology sequence of the space pair
$((\Lambda S^n)_{S^1}, \mathcal{U}_{S^1})$.
This implies there is a cocycle in $\eta$ such that $\eta(\tau)=0$
for any simplex $\tau\subset(\Lambda S^n)_{S^1}$ with support $|\tau|\subset \mathcal{U}_{S^1}$.
Thus $\gamma^{\prime}\cap\eta$
is a cycle in $\sigma_k$ with support
$|\gamma^{\prime}\cap\eta|\subset (\Lambda^{\lambda_k+\epsilon}\setminus\mathcal{U})_{S^1}$.
Now $[\widetilde{F}(1, \gamma^{\prime}\cap\eta)]=[\gamma^{\prime}\cap\eta]$
and $|\widetilde{F}(1, \gamma^{\prime}\cap\eta)|\subset\Lambda^{\lambda_k-\epsilon}$.
This contradicts to the definition of $\lambda_k$.
Hence we must have $i^\ast(\eta)\neq0$.

Now suppose
there are finitely many prime closed geodesics on $(S^n, F)$.
Then we have $E^{-1}(\lambda_k)\cap crit(E)$ is a  union of
finite many disjoint circles. Thus we can choose  $\mathcal{U}$
to be a union of finite many tubular neighborhoods
of this circles. This implies $i^\ast(\eta)=0$. This contradiction
proves the proposition.\hfill\hb

We call a closed geodesic $c$ {\it isolated}, if there exists
an $S^1$-invariant neighborhood $\mathcal{N}$ of $S^1\cdot v$ such that
$crit(E)\cap\mathcal{N}=S^1\cdot c$.

{\bf Proposition 2.3.} {\it Suppose there exists $\delta_k>0$ such that
any closed geodesic $c$ with $E(c)\in(\lambda_k-\delta_k, \lambda_k+\delta_k)$
is isolated. Then there exists a closed geodesic $c_k$
such that
$$E(c_k)=\lambda_k,\qquad i(c_k)\le n-1+2k\le i(c_k)+\nu(c_k)$$
for $0\le k\le n-1$.
 }

{\bf Proof.} By assumption and the Palais-Smale condition, we can choose
$\epsilon>0$ small enough such that $\lambda_k$ is the unique critical value of $E$
in $(\lambda_k-\epsilon, \lambda_k+\epsilon)$.

By definition of $\lambda_{k}$, we have a chain
$\gamma\in\sigma_{k}$ such that $\sup_{d\in\im\gamma}E(d)<\lambda_k+\epsilon$.
Thus we have $0\neq[\gamma]\in H_{S^1, n-1+2k}(\Lambda^{\lambda_k+\epsilon}, \Lambda^0)$.
If $H_{S^1, n-1+2k}(\Lambda^{\lambda_k+\epsilon}, \Lambda^{\lambda_k-\epsilon})=0$,
then by the exactness of the homology sequence
of the triple $((\Lambda^{\lambda_k+\epsilon})_{S^1}, (\Lambda^{\lambda_k-\epsilon})_{S^1}, (\Lambda^0)_{S^1})$,
we have $\gamma^\prime\in H_{S^1, n-1+2k}(\Lambda^{\lambda_k-\epsilon}, \Lambda^0)$.
such that $i_\ast(\gamma^\prime)=\gamma$, where
$i: ((\Lambda^{\lambda_k-\epsilon})_{S^1}, (\Lambda^0)_{S^1})\rightarrow ((\Lambda^{\lambda_k+\epsilon})_{S^1}, (\Lambda^0)_{S^1})$
is the inclusion. This contradicts to the definition of $\lambda_k$.
Hence we have
\be H_{S^1, n-1+2k}(\Lambda^{\lambda_k+\epsilon}, \Lambda^{\lambda_k-\epsilon})\neq0.
\lb{2.8}\ee
Since $\lambda_k$ is the unique critical value of $E$
in $(\lambda_k-\epsilon, \lambda_k+\epsilon)$.
By the equivariant version of Lemma 1.4.2 of \cite{Cha}, we have
\be H_{S^1, \ast}(\Lambda^{\lambda_k+\epsilon}, \Lambda^{\lambda_k-\epsilon})
\simeq\bigoplus_{c\in crit(E)} C_{S^1,\ast}(E, S^1\cdot c),
\lb{2.9}\ee
where $C_{S^1,\ast}(E, S^1\cdot c)$ is the $S^1$-critical group at $S^1\cdot c$
defined by
\bea C_{S^1,\; \ast}(E, \;S^1\cdot c)
=H_{S^1,\; \ast}(\Lambda(c)\cap\mathcal{N},\;
(\Lambda(c)\setminus S^1\cdot c)\cap\mathcal{N}),\lb{2.10}
\eea
where $\mathcal{N}$ is any $S^1$-invariant
open neighborhood of $S^1\cdot c$ such that
$crit(E)\cap\mathcal{N}=S^1\cdot c$.
By shrinking $\mathcal{N}$ if necessary, we may assume
the multiplicities $m(d)$ for $d\in \mathcal{N}$
is bounded from above. Then by Lemma 6.11 of \cite{FaR}, we have
\bea H^{\ast}_{S^1}(\Lambda(c)\cap\mathcal{N},\;
(\Lambda(c)\setminus S^1\cdot c)\cap\mathcal{N})\lb{2.10}
\simeq H^{\ast}((\Lambda(c)\cap\mathcal{N})/S^1,\;
((\Lambda(c)\setminus S^1\cdot c)\cap\mathcal{N})/S^1).
\lb{2.11}
\eea
Now by introducing finite-dimensional approximations
to $\Lambda$ and apply Gromoll-Meyer theory as in
\S6 of \cite{Rad2}, we have
\be H^{q}((\Lambda(c)\cap\mathcal{N})/S^1,\;
((\Lambda(c)\setminus S^1\cdot c)\cap\mathcal{N})/S^1)=0
\lb{2.12}
\ee
provided $q>i(c)+\nu(c)$ or $q<i(c)$, cf. Satz 6.13 of \cite{Rad2}.

Now the proposition follows from (\ref{2.8})-(\ref{2.12}).
\hfill\hb

{\bf Proposition 2.4.} {\it Suppose a closed geodesic
$c$ is non-degenerate, i.e., $\nu(c)=0$, then $c$
is isolated. }

{\bf Proof.}  Following \cite{Rad2}, Section
6.2, we introduce finite-dimensional
approximations to $\Lambda$. We choose an arbitrary energy value
$a>0$ and $k\in\N$ such that every $F$-geodesic of length
$<\sqrt{2a/k}$ is minimal. Then
$$ \Lm(k,a)=\left\{\ga\in\Lm \mid E(\ga)<a \mbox{ and }
    \ga|_{[i/k,(i+1)/k]}\mbox{ is an $F$-geodesic for }i=0,\ldots,k-1\right\} $$
is a $(k\cdot\dim M)$-dimensional submanifold of $\Lambda$
consisting of closed geodesic polygons with $k$ vertices. The set
$\Lambda(k,a)$ is invariant under the subgroup $\Z_k$ of $S^1$.
Then $E|_{\Lm(k,a)}$ is smooth.
Closed geodesics in $\Lambda^{a-}=\{\gamma\in\Lambda\mid
E(\gamma)<a\}$ are precisely the critical points of
$E|_{\Lm(k,a)}$, and for every closed geodesic $c\in\Lm(k,a)$ the
index of $(E|_{\Lm(k,a)})''(c)$ equals $i(c)$ and the null space
of $(E|_{\Lm(k,a)})''(c)$ coincides with the null space of
$E''(c)$, cf. \cite{Rad2}, p.51. Clearly,
$S^1\cdot c\subset \Lm(k,a)$ is a critical manifold
and its tangent space contains in the null space
of $(E|_{\Lm(k,a)})''(c)$. Since $\nu(c)=0$ by assumption,
we have $(E|_{\Lm(k,a)})''(c)|_{N(S^1\cdot c)}$ is non-degenerate.
where $N(S^1\cdot c)$ is the normal bundle of $S^1\cdot c$ in
$\Lm(k,a)$. Thus $S^1\cdot c$ is a non-degenerate
critical manifold of $E|_{\Lm(k,a)}$ in the sense of \cite{Bot1}.
Hence it is an isolated critical manifold of $E|_{\Lm(k,a)}$,
and then it must be a isolated critical manifold of $E$.\hfill\hb

\setcounter{equation}{0}%\setcounter{figure}{0}
\section{Proof of the main theorems}%{Section 1}

In this section, we give the proofs of the main theorems.

{\bf Proof of Theorem 1.2.} Suppose the Finsler  $n$-sphere
$(S^n,\,F)$ with reversibility $\lambda$
satisfying $F^2<(\frac{\lambda+1}{\lambda})^2g_0$
and $l(S^n, F)\ge \pi(1+\frac{1}{\lambda})$.
Then we have
$\Lambda_0^{2\pi^2+\epsilon} S^n\subset \Lambda^{2\pi^2(\frac{\lambda+1}{\lambda})^2-\epsilon} S^n$
for some $\epsilon>0$ is sufficiently small, where we use notations in
Proposition 2.1.
Thus by Proposition 2.2, there are $n$ critical values
\be 0<\lambda_0\le\lambda_1\le\cdots\le\lambda_{n-1}< 2\pi^2\left(\frac{\lambda+1}{\lambda}\right)^2.
\lb{3.1}\ee  of $E$. Hence we can find $n$ closed geodesics $\{c_k\}_{0\le k\le n-1}$
of $(S^n, F)$ with $E(c_k)=\lambda_k$. Clearly each $c_k$ is
a prime closed geodesic without self-intersections, since otherwise
we should have
\be E(c_k)\ge\frac{1}{2}L(c_k)^2\ge2\pi^2\left(\frac{\lambda+1}{\lambda}\right)^2.
\lb{3.2}\ee
Now if $0<\lambda_0<\lambda_1<\cdots<\lambda_{n-1}$,
the geodesics $\{c_k\}_{0\le k\le n-1}$ must be distinct since their energies are different.
Otherwise by Proposition 2.2, there are infinitely many non-constant closed
geodesics below the level set $E^{-1}(\lambda_{n-1}+\epsilon)$.
These proves the theorem.
\hfill\hb

{\bf Proof of Theorem 1.3.} The theorem follows directly from
Theorem 3, 4 of \cite{Rad4} and Theorem 1.2. \hfill\hb

{\bf Proof of Theorem 1.4.} Suppose the condition
$\left(\frac{\lambda}{\lambda+1}\right)^2\le \delta<K\le 1$
holds. Then Theorem 3 of \cite{Rad4} implies
$L(c_k)\ge\pi(1+\frac{1}{\lambda})$.

We have the following two cases.

{\bf Case 1.} {\it We have $0<\lambda_0<\lambda_1<\cdots<\lambda_{n-1}$.}

If there does not exist $\delta_k>0$ such that
any closed geodesic $c$ with $E(c)\in(\lambda_k-\delta_k, \lambda_k+\delta_k)$
is isolated. Then by the Palais-Smale condition and Proposition 2.4, we
can find a degenerate closed geodesic $c_k$ with $E(c_k)=\lambda_k$.
In particular, $c_k$ is non-hyperbolic.

Otherwise we can apply Proposition 2.3 to
find a closed geodesic $c_k$ satisfying
$E(c_k)=\lambda_k$ and $i(c)\le n-1+2k\le i(c)+\nu(c)$.
Clearly, in order to prove the theorem,
it is sufficient to assume $c_k$ is hyperbolic for $0\le k\le n-1$,
then we have $i(c_k)=n-1+2k$.

Choose $p, q\in\N$ such that
$\frac{p}{q}<\frac{\lambda+1}{\lambda}\sqrt{\delta}$, then we have
$L(c_k^q)\ge q\pi(1+\frac{1}{\lambda})>\frac{p\pi}{\sqrt{\delta}}$.
Hence by the  Morse-Schoenberg comparison theorem(cf. P. 220 of \cite{BTZ1}), we have
$i(c_k^q)\ge p(n-1)$.
By Corollary 2.3 of \cite{BTZ1},
we have $i(c_k^q)=q(n-1+2k)$. Thus we must have
$\frac{p}{q}\le \frac{n-1+2k}{n-1}$.
Thus if $\frac{2n-3}{n-1}\frac{\lambda}{\lambda+1}<\sqrt{\delta}\le 1$,
the geodesics $c_k$ for $0\le k< [\frac{n}{2}]$ can't be
hyperbolic.

Suppose $F^2<\frac{g_0}{\rho}$, then by the proof of Proposition2.2,
we have $L(c_k)=\sqrt{2E(c_k)}<\frac{2\pi}{\sqrt{\rho}}$. Choose $p,
q\in\N$ such that $L(c_k)<\frac{p\pi}{q}<\frac{2\pi}{\sqrt{\rho}}$.
Then $L(c_k^q)<p\pi$, and hence $i(c_k^q)+\nu(c_k^q)\le p(n-1)$
by the  Morse-Schoenberg comparison theorem. Since
$c_k$ is hyperbolic, we have $i(c_k)=n-1+2k$ and
$i(c_k^q)=q(n-1+2k)$. Hence we must have $\frac{p}{q}\ge
\frac{n-1+2k}{n-1}$. Thus if $\frac{2(n-1)}{2n-1}<\sqrt{\rho}\le 1$,
the geodesics $c_k$ for $[\frac{n+1}{2}]\le k\le n-1$ can't be
hyperbolic.

Hence the theorem holds in this case.

{\bf Case 2.} {\it We have $\lambda_k=\lambda_{k+1}$ for some $0\le k< n-1$.}

We claim that there must be infinitely many non-hyperbolic prime closed geodesics
on $(S^n, F)$.

Suppose the claim does not hold. Then there are finitely
many non-hyperbolic closed geodesics $\{d_1, \ldots,d_l\}$
in $\Lambda^{\lambda_k+1}$. Choose disjoint $S^1$-invariant
open neighborhoods $\mathcal{U}_1,\ldots,\mathcal{U}_l$ such that
$S^1\cdot d_i\subset\mathcal{U}_i$ and $\mathcal{U}_i$ is
$S^1$-homotopic to $S^1\cdot d_i$. By shrinking $\mathcal{U}_i$
if necessary, we may assume $\partial\mathcal{U}_i\cap crit(E)=\emptyset$
for each $i$. Note that there are finitely many non-constant closed geodesics
in $\Lambda^{\lambda_k+1}\setminus (\cup_{1\le i\le l}\mathcal{U}_i)$ by
Proposition 2.4. Denote them by $\{d_{l+1}, \ldots,d_m\}$.
Choose disjoint $S^1$-invariant
open neighborhoods $\mathcal{U}_{l+1},\ldots,\mathcal{U}_m$ such that
$S^1\cdot d_i\subset\mathcal{U}_i$ and $\mathcal{U}_i$ is
$S^1$-homotopic to $S^1\cdot d_i$ for $l+1\le i\le m$.
By shrinking them if necessary, we may assume
$\mathcal{U}_i\cap\mathcal{U}_j=\emptyset$ for $1\le i, j\le m$
and $i\neq j$. Denote by $\mathcal{U}=\cup_{1\le i\le m}\mathcal{U}_i$.
Then we have $crit(E)\cap E^{-1}(\lambda_k)\subset \mathcal{U}$.
Let $i: \mathcal{U}\rightarrow \Lambda S^n$
be the inclusion and $\eta\in H^2_{S^1}(\Lambda S^n)$
as in the proof of Proposition 2.2.
Then we have $i^\ast(\eta)\neq 0$ as in Proposition 2.2.
But $\mathcal{U}$ is $S^1$-homotopic to a union of
a finite number of disjoint circles.
Hence we have $i^\ast(\eta)=0$.
This contradiction proves
there must be infinitely many non-hyperbolic prime closed geodesics
on $(S^n, F)$.

Combining the above two cases, we obtain the theorem.\hfill\hb

\bibliographystyle{abbrv}

\bigskip

\end{document}